\documentclass[review]{elsarticle}
\usepackage{amssymb}
\usepackage{amsthm}
\usepackage{amsmath}
\usepackage{graphicx}
\usepackage{comment}
\usepackage[usenames,dvipsnames]{xcolor}
\usepackage{indentfirst,latexsym,bm,amsmath,amssymb,amsthm}
\usepackage{helvet} %helvetica
\newtheorem{thm}{Theorem}[section]

\theoremstyle{definition}

\theoremstyle{remark}
\newtheorem*{rem}{Remark}

\numberwithin{equation}{section}

\theoremstyle{remark}

\makeatletter

\newcommand{\Rmnum}[1]{\expandafter\@slowromancap\romannumeral #1@}
\makeatother
\newcommand{\grad}{\nabla}   %ÌÝ¶È
\newcommand{\divg}{\nabla\cdot}  %É¢¶È
\begin{document}

\begin{frontmatter}

\title{Global Smooth Solutions of MHD equations with large data}

\author{Yurui Lin$^1$}
%\address{
	%Cangshan District, 	Fuzhou}
\fntext[]{College of Computer and Information Science, Fujian Agriculture and Forestry University, China.}

\author{Huali Zhang$^2$}
%\address{220 Handan Rd., Yangpu District, Shanghai}
\fntext[myfootnote]{School of Mathematical Sciences, Fudan University,China}

\author{Yi Zhou$^{3*}$}
%\address{220 Handan Rd., Yangpu District, Shanghai}
\fntext[]{$^*$Corresponding author: {\tt
		Email address:yizhou@fudan.edu.cn}, College of Computer and Information Science, Fujian Agriculture and Forestry University, China and School of Mathematical sciences, Fudan Universety, China}
%[mysecondaryaddress]{School of Mathematical sciences, Fudan Universety, China\corref{mycorrespondingauthor}}
%\ead{yizhou@fudan.edu.cn}

\begin{abstract}
In this paper, we establish the global existence of smooth solutions of the three-dimensional MHD system for a class of large initial data. Both the initial velocity and magnetic field can be arbitrarily large in the critical norm.
\end{abstract}


\end{frontmatter}

%\linenumbers

\section{Introduction}
\qquad
Magnetohydrodynamics (MHD) is to study the dynamics of the velocity and magnetic fields in electrically conducting fluids such as  plasmas, liquid metals, and salt water or electrolytes.

\qquad
In this paper, we consider the Cauchy's problem for the incompressible MHD system in three space dimensions
\begin{align}\label{1.1}
	&u_t+u\cdot\nabla u-h\cdot\nabla h+\nabla q=\nu\Delta u \ \ \ \ \ \ t>0,\ x\in \mathbb{R}^3,\\
	&h_t+u\cdot\nabla h-h\cdot\nabla u=\mu\Delta h \ \ \ \ \ \ \ \ \ \ \ \ \ t>0,\ x\in \mathbb{R}^3, \\
	&\nabla\cdot u=\nabla\cdot h=0 \ \ \ \ \ \ \ \ \ \ \ \ \ \ \ \ \ \ \ \ \ \ \ \ \ \ \ t>0,\ x\in \mathbb{R}^3, \\
	&t=0: \quad u=u_0(x), \quad h=h_0(x), \ \ \ \ \ \ \ \ \ x\in \mathbb{R}^3,
\end{align}
where
$u$ denotes the fluid velocity, $h$ the magnetic field, $q$ the pressure, $\nu$ the viscosity, $\mu$ the magnetic diffusivity, $\nu$ and $\mu$ are positive constants. $u_0$ and $h_0$ are the initial data which satisfy
\begin{equation}\label{1.5}
	\nabla\cdot u_0 =\nabla\cdot h_0=0, \qquad x\in \mathbb{R}^3,
\end{equation}
and
\begin{equation}\label{1.6}
	u_0,\ h_0 \in C^{\infty}(\mathbb{R}^3).
\end{equation}
The equation (1.1) is the law of conservation of momentum, (1.2) is the induction equation and (1.3) reflects the conservation of mass. Besides its widely applications in physics,
it remains an outstanding  mathematical  problem whether  there always exists a global smooth solution for smooth initial data.

\qquad In the case that the magnetic field $h$ is a constant vector, that is, in the case of the Navier-Stokes equations, the global well-posedness of strong solutions in three space dimensions has been studied by many authors for small initial data in various spaces, for example, in Sobolev space ${\rm \dot{H}}^s$, $s \geq \frac{1}{2}$ by Fujita and Kato \cite{FK} (see also Chemin \cite{Chemin}), in Lebesgue space $L^3(\mathbb{R}^3)$ by Kato \cite{Kato}, in Besov spaces $\dot{B}^{- 1 + \frac{3}{p}}_{p, \infty}$, $n < p < \infty$ by Cannone \cite{Cannone} and Planchon \cite{Planchon}, and in the larger space ${\rm BMO}^{- 1}$ by Koch and Tataru \cite{KT}. These small data global existence results can
be extended to MHD system with non-trivial magnetic field $h$ provided that both the initial velocity field $u_0$ and the initial
magnetic field $h_0$ satisfy appropriate smallness conditions. For example, see \cite{Duv, MYZ, Sermange}.

\qquad
The MHD system (1.1)-(1.4) is called the ideal MHD system if $\nu=0$ and $\mu>0$. In \cite{Lei1}, Lei proved the global regularity of axially symmetric solutions to the system of incompressible ideal magnetohydrodynamics in three space dimensions by inspiration of the axially symmetric Navier Stokes equations, while the initial data is not required to be small. The MHD system (1.1)-(1.4) is called non-resistive MHD system if $\mu=0$ and $\nu > 0$. Lin and Zhang \cite{Lin2} established the global well-posedness of the imcompressible MHD type system with initial data close to some equilibrium states in three space dimensions. Similar results were also obtained in two space dimensions by Lin et al \cite{Lin1}. Zhang \cite{Zhangting}, Ren et al \cite{Ren} provided two simplified proofs of the main result in \cite{Lin1}. For the MHD system (1.1)-(1.4) with $\mu=0$ and $\nu =0$, Bardos et al \cite{Bar} proved the global existence of the smooth solution when the initial data is close to the equilibtium state.

\qquad
The aim of this paper is to construct global smooth solutions to the Cauchy's problem of MHD system (1.1)-(1.4) with a class of large initial data in three space dimensions. Our result generalizes that of Lei et.al \cite{Lei3} for the 3-D incompressible Navier-Stokes equations to the case of MHD with large velocity fields and large magnetic fields. We also give a simple and completely different proof than that of Lei et al \cite{Lei3}. Their proof relies on energy estimates and works in the physical space variables while our proof works in the Fourier Transforms variables. We first establish a stability result for the MHD system which generalizes a small data global well-posedness result of Lei and Lin \cite{Lei2} for the 3-D incompressible Navier-Stokes equations. We then explore the structure of the nonlinear term of the MHD to show that our data is 'nonlinearly small' in the sense of Chemin and Gallagher\cite{Chemin1}. Our main result can be stated as follows.

\begin{thm}
	Consider the Cauchy problem (1.1)-(1.4). Suppose that
	\begin{equation}\label{1.7}
		u_0(x)=u_{01}(x)+u_{02}(x),
	\end{equation}
	\begin{equation}\label{1.8}
		h_0(x)=h_{01}(x)+h_{02}(x),
	\end{equation}
	with
	\begin{equation}\label{1.111}
		\divg u_{01}(x)=\divg h_{01}(x)=0,
	\end{equation}
	\begin{equation}\label{1.9}
		u_{02}(x)=\alpha_1 v_0(x),
	\end{equation}
	\begin{equation}\label{1.10}
		h_{02}(x)=\alpha_2 v_0(x),
	\end{equation}
	where $\alpha_1$, $\alpha_2$ are two real constants and $v_0(x)$ has the following properties

	\qquad $ (1) \ \divg v_0(x)=0,$

	\qquad$ (2) \ \grad\times v_0(x) =\sqrt{-\Delta} v_0(x),$

	\qquad$ (3) \ supp \widehat{v_0}(\xi)\subset \{ \xi : 1-\delta \leq |\xi|\leq 1+\delta,  \},  \ 0\leq\delta\leq \frac{1}{2},$
	
   \qquad $(4)\ |v_0|_{\mathbf{L}^2(\mathbb{R}^3 )} \leq M\delta^{-\frac{1}{2}}.$

where $\widehat{v_0}$ denote the Fourier transform of $v_0$.

	\qquad Then there exists an positive constant $\varepsilon$ sufficiently small, and depending only on $\mu$, $\nu$, $k$, $\alpha_1$, $\alpha_2$ and $M$ such that the Cauchy problem (1.1)-(1.4) has a global smooth solution provided that
	\begin{equation}\label{1.11}
		\delta\leq \varepsilon
	\end{equation}
	and
	\begin{equation}\label{1.12}
		\int_{\mathbb{R}^3} \frac{|\widehat{u_{01}}(\xi)|}{|\xi|}d\xi +\int_{\mathbb{R}^3} \frac{|\widehat{h_{01}}(\xi)|}{|\xi|}d\xi \leq \varepsilon.
	\end{equation}
\end{thm}

\begin{rem}
	$v_0(x)$ satisfies $(1)-(4)$ have been constructed by Lei et al \cite{Lei2}. In the limiting case of $\delta=0$, $\grad \times v_0=v_0$,
	and the flow is called Beltrami flow.
\end{rem}
\begin{rem}
	The constant $M$ can be arbitrary large, so our initial data can be arbitrary large .
\end{rem}

Finally, we remark that it is an open problem if there would still have global existence for the type of initial data in Theorem 1.1, if $\alpha_2=1$, $\mu=0$.

\section{Proof of the main result}
By Schwartz' inequality,
\begin{equation}
\ |\hat{v}_0|_{\mathbf{L}^1(\mathbb{R}^3 )}\le C\delta^{\frac{1}{2}}\ |\hat{v}_0|_{\mathbf{L}^2(\mathbb{R}^3 )}=C\delta^{\frac{1}{2}}\ |v_0|_{\mathbf{L}^2(\mathbb{R}^3 )}\le CM
\end{equation}

\qquad Let $f$ satisfies the heat equation
\begin{equation}\label{2.1}
	\begin{cases}
		f_t-\nu \Delta f=0\\
		t=0:\quad f=u_{02}(x),
	\end{cases}
\end{equation}
and the condition
$$ \divg u_{02}=0 $$
implies
\begin{equation}\label{2.2}
	\divg f=0.
\end{equation}
\qquad In a similar way, let $g$ satisfies the heat equation
\begin{equation}\label{2.3}
	\begin{cases}
		g_t-\mu \Delta g=0\\
		t=0:\quad g=h_{02}(x),
	\end{cases}
\end{equation}
and the condition
$$ \divg h_{02}=0$$
implies
\begin{equation}\label{2.4}
	\divg g=0.
\end{equation}
\qquad In the beginning of this section, we do not require $u_{02}$, $h_{02}$ be the the form of \eqref{1.9} \eqref{1.10}, we only assume that
\begin{equation}\label{2.5}
	\int_{\mathbb{R}^3} (1+\frac{1}{|\xi|^2})|\widehat{u_{02}}(\xi)|d\xi+\int_{\mathbb{R}^3} (1+\frac{1}{|\xi|^2})|\widehat{h_{02}}(\xi)|d\xi \leq 4CM,
\end{equation}
let
$$ U=u-f, \quad H=h-g, $$
then (1.1)(1.2) can be written as follows
\begin{equation}\label{2.6}
	U_t+U\cdot\nabla U+f\cdot\nabla U+U\cdot\nabla f-H\cdot\nabla H-g\cdot\nabla H-H\cdot\nabla g+\nabla q= \nu \Delta U+F_1,
\end{equation}
\begin{equation}\label{2.7}
	H_t+U\cdot\nabla H+f\cdot\nabla H+U\cdot\nabla g-H\cdot\nabla U-g\cdot\nabla U-H\cdot\nabla f= \mu \Delta H+G,
\end{equation}
where
\begin{equation}\label{2.8}
	\begin{split}
		F_1=& g\cdot\nabla g-f\cdot\nabla f\\
		=&g\times (\nabla\times g)-f\times (\nabla\times f)+\nabla (\frac{1}{2}|f|^2-\frac{1}{2}|g|^2)\\
		\doteq& F+\nabla (\frac{1}{2}|f|^2-\frac{1}{2}|g|^2),
	\end{split}
\end{equation}
\begin{equation}\label{2.90}
	F=g\times (\nabla\times g)-f\times (\nabla\times f),
\end{equation}
\begin{equation}\label{2.9}
	G= g\cdot\nabla f-f\cdot\nabla g = \nabla\times (g\times f),
\end{equation}
absorbing $\nabla (\frac{1}{2}|f|^2-\frac{1}{2}|g|^2)$ into the presure term and noting the divergence free condition, we finally get
\begin{equation}\label{2.10}
	U_t-\nu \Delta U=-\mathbb{P}\nabla \cdot(U\otimes U+f\otimes U+U\otimes f-H\otimes H-g\otimes H-H\otimes g)+\mathbb{P}F,
\end{equation}
where
$$ \mathbb{P}=(I-(-\Delta)^{-1}\nabla\times\nabla\times)$$
is the projection to the divergence free fields. Similarly,
\begin{equation}\label{2.11}
	H_t-\mu\Delta H=-\nabla\cdot(U\otimes H+f\otimes H+U\otimes g-H\otimes U-g\otimes U-H\otimes f)+G
\end{equation}
\begin{equation}\label{2.12}
	t=0: \quad U=u_{01}(x), \quad H=h_{01}(x)
\end{equation}
\begin{equation}\label{2.13}
	\nabla\cdot U=\nabla\cdot H=0
\end{equation}
we first prove a stability result which generalizes a theorem of Lei and Lin \cite{Lei2}.
\begin{thm}
	Consider the system of equations (\ref{2.10})-(\ref{2.13}), suppose that $f$ satisfies (\ref{2.1})(\ref{2.2}) and $g$ satisfies (\ref{2.3})(\ref{2.4}) and suppose that (\ref{2.5}) holds. Then there exists a small positive constant $\delta_0$ depending only on $\mu$, $\nu$ and $M$ such that (\ref{2.10})-(\ref{2.13}) has a global smooth solution satisfies
	\begin{equation}\label{2.14}
		\begin{split}
			&\int_{\mathbb{R}^3}\frac{|\widehat{U}(t,\xi)|}{|\xi|}d\xi+\nu\int_0^t \int_{\mathbb{R}^3}|\xi||\widehat{U}(\tau,\xi)|d\xi d\tau \\
			+&\int_{\mathbb{R}^3}\frac{|\widehat{H}(t,\xi)|}{|\xi|}d\xi+\mu \int_0^t \int_{\mathbb{R}^3}|\xi||\widehat{H}(\tau,\xi)|d\xi d\tau \leq C_{*}\delta_0,
		\end{split}
	\end{equation}
	provided that
	\begin{equation}\label{2.15}
		\begin{split}
			 &\int_{\mathbb{R}^3}\frac{|\widehat{u_{01}}(\xi)|}{|\xi|}d\xi+\int_{\mathbb{R}^3}\frac{|\widehat{h_{01}}(\xi)|}{|\xi|}d\xi\\
			+&\int_0^t \int_{\mathbb{R}^3}\frac{|\widehat{F}(\tau,\xi)|}{|\xi|}d\xi d\tau+\int_0^t \int_{\mathbb{R}^3}\frac{|\widehat{G}(\tau,\xi)|}{|\xi|}d\xi d\tau \leq \delta_0.
		\end{split}
	\end{equation}
\end{thm}

\qquad \textit{Proof.} Let
\begin{equation}\label{2.16}
	 E_0(t)=\int_{\mathbb{R}^3}\frac{|\widehat{U}(t,\xi)|}{|\xi|}d\xi+\int_{\mathbb{R}^3}\frac{|\widehat{H}(t,\xi)|}{|\xi|}d\xi
	,
\end{equation}
\begin{equation}\label{2.17}
	E_1(t)=\nu\int_0^t \int_{\mathbb{R}^3}|\xi||\widehat{U}(\tau,\xi)|d\xi d\tau +\mu \int_0^t \int_{\mathbb{R}^3}|\xi||\widehat{H}(\tau,\xi)|d\xi d\tau,
\end{equation}
taking Fourier transform of (\ref{2.10}) (\ref{2.11}) we get
\begin{equation}\label{311}
	\partial_t\hat{U}+\nu |\xi|^2\hat{U}=\hat{A}+\hat{\mathbb{P}F}
\end{equation}
and
\begin{equation}\label{322}
	\partial_t\hat{H}+\mu |\xi|^2\hat{H}=\hat{B}+\hat{G},
\end{equation}
where
\begin{equation*}
	A=-\mathbb{P}\nabla \cdot(U\otimes U+f\otimes U+U\otimes f-H\otimes H-g\otimes H-H\otimes g),
\end{equation*}
\begin{equation*}
	B=-\nabla\cdot(U\otimes H+f\otimes H+U\otimes g-H\otimes U-g\otimes U-H\otimes f).
\end{equation*}
taking a inner product of (\ref{311}), (\ref{322}) by $\frac{\bar{\hat{U}}}{|\hat{U}|}\frac{1}{|\xi|}$, $\frac{\bar{\hat{H}}}{|\hat{H}|}\frac{1}{|\xi|}$ respectively, we get
\begin{equation}\label{4111}
	\frac{1}{|\xi|}\partial_t\hat{U}\cdotp \frac{\bar{\hat{U}}}{|\hat{U}|}+ \nu |\xi||\hat{U}|=\frac{\hat{A}}{|\xi|} \cdot
	\frac{\bar{\hat{U}}}{|\hat{U}|}+\frac{\hat{\mathbb{P}F}}{|\xi|} \cdot
	\frac{\bar{\hat{U}}}{|\hat{U}|},
\end{equation}
and
\begin{equation}\label{4222}
	\frac{1}{|\xi|}\partial_t\hat{H}\cdotp \frac{\bar{\hat{H}}}{|\hat{H}|}+\mu |\xi||\hat{H}|=\frac{\hat{B}}{|\xi|} \cdot \frac{\bar{\hat{H}}}{|\hat{H}|}+\frac{\hat{G}}{|\xi|} \cdot \frac{\bar{\hat{H}}}{|\hat{H}|},
\end{equation}
noting $Re(\partial_t\hat{U}\cdotp \frac{\bar{\hat{U}}}{|\hat{U}|})=\partial_t(|\hat{U}|)$, $Re(\partial_t\hat{H}\cdotp \frac{\bar{\hat{H}}}{|\hat{H}|})=\partial_t(|\hat{H}|)$, taking the real part of \eqref{4111} and \eqref{4222} we obtain
\begin{equation}\label{411}
	\partial_t(\frac{|\hat{U}|}{|\xi|})+\nu |\xi||\hat{U}|=Re(\frac{\hat{A}}{|\xi|} \cdot
	\frac{\bar{\hat{U}}}{|\hat{U}|}+\frac{\hat{\mathbb{P}F}}{|\xi|} \cdot
	\frac{\bar{\hat{U}}}{|\hat{U}|})
\end{equation}
and
\begin{equation}\label{422}
	\partial_t({\frac{|\hat{H}|}{|\xi|}})+\mu |\xi||\hat{H}|=Re(\frac{\hat{B}}{|\xi|} \cdot \frac{\bar{\hat{H}}}{|\hat{H}|}+\frac{\hat{G}}{|\xi|} \cdot \frac{\bar{\hat{H}}}{|\hat{H}|}),
\end{equation}
obviously,
\begin{equation*}
	\frac{|\hat{\mathbb{P}F}|}{|\xi|}\leq \frac{|\hat{F}|}{|\xi|},
\end{equation*}
\begin{equation*}
	\frac{|\hat{A}|}{|\xi|}\leq |\hat{U}\ast \hat{U}|+2|\hat{f}\ast \hat{U}|+|\hat{H}\ast \hat{H}|+2|\hat{H}\ast \hat{g}|,
\end{equation*}
\begin{equation*}
	\frac{|\hat{B}|}{|\xi|}\leq |\hat{U}\ast \hat{H}|+|\hat{f}\ast \hat{H}|+|\hat{U}\ast \hat{g}|+|\hat{H}\ast \hat{U}|+|\hat{U}\ast \hat{g}|+|\hat{H}\ast \hat{f}|.
\end{equation*}
Adding (\ref{411}) (\ref{422}) and integrating on $[0, t]\times \mathbb{R}^3$, noting (\ref{2.15}) we obtain
\begin{equation}\label{2.18}
	\begin{split}
		E_0(t)+E_1(t) \leq \delta_0+&C\int_0^t\int_{\mathbb{R}^3}\int_{\mathbb{R}^3} (|\widehat{U}(\tau,\xi-\eta)|+|\widehat{H}(\tau,\xi-\eta)|)\cdot (|\widehat{U}(\tau,\eta)|+|\widehat{H}(\tau,\eta)|)d\xi d\eta d\tau \\
		 +&C\int_0^t\int_{\mathbb{R}^3}\int_{\mathbb{R}^3}(|\widehat{U}(\tau,\xi-\eta)|+|\widehat{H}(\tau,\xi-\eta)|)(|\widehat{f}(\tau,\eta)|+|\widehat{g}(\tau,\eta)|)d\xi d\eta d\tau \\
		=& \delta_0+\mathrm{\Rmnum{1}}+\mathrm{\Rmnum{2}}.
	\end{split}
\end{equation}
Here
\begin{equation*}
	\mathrm{\Rmnum{1}}=C\int_0^t\int_{\mathbb{R}^3}\int_{\mathbb{R}^3} (|\widehat{U}(\tau,\xi-\eta)|+|\widehat{H}(\tau,\xi-\eta)|)\cdot (|\widehat{U}(\tau,\eta)|+|\widehat{H}(\tau,\eta)|)d\xi d\eta d\tau,
\end{equation*}
\begin{equation*}
	 \mathrm{\Rmnum{2}}=C\int_0^t\int_{\mathbb{R}^3}\int_{\mathbb{R}^3}(|\widehat{U}(\tau,\xi-\eta)|+|\widehat{H}(\tau,\xi-\eta)|)(|\widehat{f}(\tau,\eta)|+|\widehat{g}(\tau,\eta)|)d\xi d\eta d\tau.
\end{equation*}
By using the technics of Lei and Lin \cite{Lei2}
\begin{equation}\label{2.19}
	2 \leq |\xi-\eta||\eta|^{-1}+|\xi-\eta|^{-1}|\eta|,
\end{equation}
we get
\begin{equation}\label{2.20}
	I \leq C \sup_t E_0(t)E_1(t).
\end{equation}
\qquad To estimate \Rmnum{2}, we divide the case $|\xi-\eta|\leq L$ and $|\xi-\eta|\geq L$, then
\begin{equation}\label{2.21}
	\begin{split}
		\mathrm{\Rmnum{2}}\leq & CL\int_0^t \big{(} \int_{\mathbb{R}^3} (|\widehat{f}(\tau,\eta)|+|\widehat{g}(\tau,\eta)| )d\eta \big{)} E_0(\tau)d\tau\\
		+&CL^{-1}\big{(}\sup_{ \tau} \int_{\mathbb{R}^3} (|\widehat{f}(\tau,\eta)|+|\widehat{g}(\tau,\eta)| )d\eta \big{)} E_1(t).
	\end{split}
\end{equation}
\qquad By (2.5), we have
\begin{equation}\label{2.22}
	\sup_{\tau} \int_{\mathbb{R}^3} (|\widehat{f}(\tau,\eta)|+|\widehat{g}(\tau,\eta)| )d\eta \leq 4M,
\end{equation}
so we get
\begin{equation}\label{2.23}
	\begin{split}
		E_0(t)+E_1(t) \leq& \delta_0+C (\sup_t E_0(t))E_1(t) \\
		&+CML^{-1} E_1(t)+CL \int_0^t \big{(}\int_{\mathbb{R}^3}(|\widehat{f}(\tau,\eta)|+|\widehat{g}(\tau,\eta)| )d\eta \big{)}E_0(\tau)d\tau.
	\end{split}
\end{equation}
We prove \eqref{2.14} by induction, we assume that
\begin{equation}\label{2.24}
	E_0(t)+E_1(t) \leq 2C_{*}\delta_0,
\end{equation}
then
\begin{equation}\label{2.25}
	C (\sup_t E_0(t))E_1(t) \leq 4CC_{*}^2\delta_0^2 \leq \delta_0,
\end{equation}
provided that $\delta_0$ is sufficiently small.

\qquad We take $L$ in such a way that
\begin{equation}\label{2.26}
	CML^{-1}=\frac{1}{2} \quad \Rightarrow L=2CM,
\end{equation}
then we finally get
\begin{equation}\label{2.27}
	E_0(t)+\frac{1}{2}E_1(t) \leq 2\delta_0+CM\int_0^t \big{(}\int_{\mathbb{R}^3}(|\widehat{f}(\tau,\eta)|+|\widehat{g}(\tau,\eta)| )d\eta \big{)}E_0(\tau)d\tau.
\end{equation}
\qquad By Gronwall's inequality, we get
\begin{equation}\label{2.28}
	\begin{split}
		&E_0(t)+E_1(t) \\
		\leq & 8\delta_0\exp \big{(}CM\int_0^t \int_{\mathbb{R}^3}(|\widehat{f}(\tau,\eta)|+|\widehat{g}(\tau,\eta)| )d\eta d\tau \big{)}\\
		\leq& 8\delta_0 \exp (CM^2),
		%\leq& C_*\delta_0,
	\end{split}
\end{equation}
take
\begin{equation}\label{2.29}
	C_{*} =8\exp(CM^2),
\end{equation}
\qquad We conclude the proof of Theorem 2.1. \hfill $\square$

\qquad Now we prove Theorem 1.1. By Theorem 2.1, we only need to estimate
$$ \int_0^{\infty}\int_{\mathbb{R}^3} \frac{|\widehat{F}(t,\xi)|}{|\xi|}d\xi dt +\int_0^{\infty}\int_{\mathbb{R}^3} \frac{|\widehat{G}(t,\xi)|}{|\xi|}d\xi dt, $$
we have
\begin{equation}\label{2.30}
	\begin{split}
		F=& g\times (\nabla\times g)-f\times (\nabla\times f)\\
		=&g\times (\sqrt{-\Delta}g)-f\times (\sqrt{-\Delta}f)\\
		=&g\times (\sqrt{-\Delta}g-g)-f\times (\sqrt{-\Delta}f-f),
	\end{split}
\end{equation}
where $\delta$ is sufficiently small, $\widehat{F}(t,\xi)$ will be supported on $|\xi|\le 4$.
\\Thus
\begin{equation}\nonumber
	\begin{split}
		&\int_0^{\infty}\int_{\mathbb{R}^3} \frac{|\widehat{F}(t,\xi)|}{|\xi|}d\xi dt \\
		\leq& C\int_0^{\infty}\big{(}\int_{\mathbb{R}^3}|\widehat{F}(t,\xi)|^2d\xi\big{)}^{1/2} dt\\
		\leq& C\int_0^{\infty} \big{(} |\widehat{g}(t)|_{\mathbf{L}^2}|(|\xi|-1)\widehat{g}(t)|_{\mathbf{L}^1}
+|\widehat{f}(t)|_{\mathbf{L}^2}|(|\xi|-1)\widehat{f}(t)|_{\mathbf{L}^1} \big{)}dt\\
		\leq&  C\delta \big{[} \big{(} \int_0^{\infty} |\widehat{g}(t)|_{\mathbf{L}^2}dt\big{)} \sup_t |\widehat{g}(t)|_{\mathbf{L}^1} +\big{(} \int_0^{\infty} |\widehat{f}(t)|_{\mathbf{L}^2}dt\big{)} \sup_t |\widehat{f}(t)|_{\mathbf{L}^1} \big{]}\\
		\leq& CM^2\delta^{\frac{1}{2}}.
	\end{split}
\end{equation}
\qquad On the other hand
\begin{equation}\nonumber
	\begin{split}
		&\int_0^{\infty}\int_{\mathbb{R}^3} \frac{|\widehat{G}(t,\xi)|}{|\xi|}d\xi dt \\
		\leq& C\int_0^{\infty}\int_{\mathbb{R}^3} |\widehat{f\times g}(t,\xi)|d\xi dt,
	\end{split}
\end{equation}
if $\mu=\nu$ then $f\times g=0$. So we consider the case $\mu<\nu$, the case $\mu>\nu$ can be considered in a similar way.
\begin{equation}\nonumber
	\begin{split}
		&\widehat{f\times g}(t,\xi)\\
		=& \alpha_1\alpha_2 \int_{\mathbb{R}^3} e^{-\nu|\xi-\eta|^2t-\mu|\eta|^2t} \widehat{v_0}(\xi-\eta)\times \widehat{v_0}(\eta) d\eta\\
		=& \alpha_1\alpha_2 \int_{\mathbb{R}^3} e^{-\nu|\eta|^2t-\mu|\xi-\eta|^2t}\widehat{v_0}(\eta)\times \widehat{v_0}(\xi-\eta) d\eta\\
		=& \frac{1}{2}\alpha_1\alpha_2 \int_{\mathbb{R}^3}\big{(}e^{-\nu|\xi-\eta|^2t-\mu|\eta|^2t} -e^{-\nu|\eta|^2t-\mu|\xi-\eta|^2t}\big{)}\widehat{v_0}(\xi-\eta)\times \widehat{v_0}(\eta) d\eta,
	\end{split}
\end{equation}
we have
\begin{equation}\nonumber
	\begin{split}
		& |e^{-\nu|\xi-\eta|^2t-\mu|\eta|^2t} -e^{-\nu|\eta|^2t-\mu|\xi-\eta|^2t}|\\
		=& e^{-\mu(|\xi-\eta|^2+|\eta|^2)t}|e^{-(\nu-\mu)|\xi-\eta|^2t}-e^{-(\nu-\mu)|\eta|^2t} |\\
		\leq & C e^{-\mu(|\xi-\eta|^2+|\eta|^2)t} ||\xi-\eta|^2-|\eta|^2|t\\
		\leq & C e^{-\frac{\mu}{2}(|\xi-\eta|^2+|\eta|^2)t}\frac{||\xi-\eta|^2-|\eta|^2|}{|\xi-\eta|^2+|\eta|^2}.
	\end{split}
\end{equation}
\qquad In the support of $\widehat{v_0}(\xi-\eta)\times \widehat{v_0}(\eta)$, we have
\begin{equation}\nonumber
	\frac{||\xi-\eta|^2-|\eta|^2|}{|\xi-\eta|^2+|\eta|^2} \leq 10\delta,
\end{equation}
thus, it is not difficult to arrive at
\begin{equation}\nonumber
	\int_0^{\infty}\int_{\mathbb{R}^3}  |\widehat{f\times g}(t,\xi)|d\xi dt \leq CM^2\delta.
\end{equation}
\qquad This concludes the proof of Theorem 1.1.\hfill $\square$

\section*{\textbf{Acknowledgement}}
\indent
The third  author is supported by Key Laboratory of Mathematics for Nonlinear Sciences
(Fudan University), Ministry of Education of China, P. R. China, Shanghai Key Laboratory
for Contemporary Applied Mathematics, School of Mathematical Sciences, Fudan University,
P. R. China, NSFC(grants No.11421061), 973 Program(grant No. 2013CB834100) and 111 project.


\section*{References}

\begin{thebibliography}{4}
	\addcontentsline{toc}{section}{References}
	%AAAAA
	
	%BBBBBB
	\bibitem{Bar}
	Bardos, C.; Sulem, C.; Sulem, P.-L. \textit{Longtime dynamics of a conductive fluid in the presence of a strong magnetic field}. Trans. Amer. Math. Soc. 305 (1988), no. 1, 175-191.
	
	%CCCCC
	\bibitem{Cannone} M. Cannone,  \textit{Ondelettes, paraproduits et Navier-Stokes}.
	Diderot editeur, Arts et Sciences, 1995.
	
	\bibitem{Chemin} J. M. Chemin,  \textit{Remarques sur l'sexistence globale pour le syst${\rm \grave{e}}$me
		de Navier-Stokes incompressible}. SIAM Journal on Mathematical
	Analysis 1992; 23:20-28.
	\bibitem{Chemin1}
	Chemin, Jean-Yves; Gallagher, Isabelle, \textit{Well-posedness and stability results for the Navier-Stokes equations in $R^3$}, Ann. Inst. H. Poincare Anal. Non Lineaire 26 (2009), no. 2, 599-624.

	%DDDDDD
	\bibitem{Duv}
	G. Duvaut and J. L. Lions, \textit{In$\mathrm{\acute{e}}$quations en thermo$\mathrm{\acute{e}}$lasticit$\mathrm{\acute{e}}$ etmagn$\mathrm{\acute{e}}$tohydrodynamique}, Arch. Ration. Mech. Anal., 46 (1972), pp. 241-279.
	(1972), 241--279.
	%FFFFFF

	\bibitem{FK} H. Fujita and T. Kato, \textit{On the Navier-Stokes initial value problem I}.
	Archive for Rational Mechanics and Analysis 1964; 16:269--315.
	%GGGGG
	
	%KKKKK
	\bibitem{Kato} T. Kato,  \textit{Strong $L^q$ solutions of the Navier-Stokes
		equations in $R^n$ with applications to weak solutions}, Mathematische
	Zeitschrift 1984; 187:471--480.
	\bibitem{KT} H. Koch and D. Tataru,  \textit{Well-posedness for the Navier-Stokes
		equations}, Advance Mathematics 2001; 157:22--35.
	%LLLLLL
	\bibitem{Lei1}
	Zhen, Lei, \textit{On Axially Symmetric Incompressible Magnetohydrodynamics in Three Dimensions},
	arXiv:1212.5968.
	\bibitem{Lei2}
	Lei, Zhen, Lin, Fanghua, \textit{Global mild solutions of Navier-Stokes equations}, Comm. Pure Appl. Math. 64 (2011), no. 9, 1297-1304.
	\bibitem{Lei3}
	Lei, Zhen, Lin, Fanghua, Zhou, Yi, \textit{Structure of Helicity and Global Solutions of Incompressible Navier-Stokes Equation}, Arch. Ration. Mech. Anal. 218 (2015), no. 3, 1417-1430..
	
	\bibitem{Lin1}
	F.H. Lin, L. Xu and P. Zhang, \textit{Global small solutions to 2-D incompressible MHD system}, arXiv:1302.5877.
	\bibitem{Lin2}
	Lin, Fanghua; Zhang, Ping,\textit{Global small solutions to an MHD-type system: the three-dimensional case}, Comm. Pure Appl. Math. 67 (2014), no. 4, 531-580.
	%PPPPPPP
	\bibitem{Planchon} F. Planchon,  \textit{Global strong solutions in Sobolev or Lebesgue spaces
		to the incompressible Navier-Stokes equations in $\mathbb{R}^3$},
	Ann. Inst. H. Poincare Anal. Non Lineaire 13 (1996), no. 3, 319-336
	%MMMMMM
	\bibitem{MYZ} C. Miao, B. Yuan and B. Zhang, \textit{Well-posedness for the
		incompressible magneto-hydrodynamic system}, Math. Methods Appl. Sci. 30 (2007), no. 8, 961--976.
	%\bibitem{Miao2}
	%C. Miao, B. Yuan, B. Zhang, Well-posedness for the incompressible magneto-hydrodynamic system, Math. Methods Appl.Sci. 30 (2007) 961-976.
	
	%RRRRRR
	\bibitem{Ren}
	Xiaoxia, Ren; Jiahong, Wu; Zhaoyin, Xiang; Zhifei, Zhang, \textit{Global existence and decay of smooth solution for the 2-D MHD equations without magnetic diffusion}. J. Funct. Anal. 267 (2014), no. 2, 503-541.
	%SSSSSS

	\bibitem{Sermange}
	M. Sermange and R. Temam, \textit{Some mathematical questions related to the MHD equations},
	Commun. Pure Appl. Math., 36 (1983), pp. 635-664.
	%XXXXXXX
	%\bibitem{Xu}
	%Xu, Li; Zhang, Ping Global small solutions to three-dimensional incompressible magnetohydrodynamical system. SIAM J. Math. Anal. 47 (2015), no. 1, 26?65.
	%ZZZZZZZ
	\bibitem{Zhangting}
	Zhang,Ting, \textit{An elementary proof of the global existence and uniqueness theorem to 2-D incompressible non-resistive MHD system}, arXiv: 1404.5681.
\end{thebibliography}
\end{document}